\pdfoutput=1
\documentclass[runningheads]{llncs}
\usepackage[T1]{fontenc}
\usepackage{graphicx}
\usepackage{caption}  
\begin{document}

\title{Equilibrium Strategies of Carbon Emission Reduction in Agricultural Product Supply Chain under Carbon Sink Trading
\thanks{This work is supported by NSFC (No.72271109).
\\$*$ Corresponding author.}}

%
%
\author{Tingting Meng\inst{1}, 
Yukun Cheng\inst{1*}, Xujin Pu\inst{1}, and
Rui Li\inst{2}
}
\authorrunning{Tingting Meng et al.}
%

\institute{
Jiangnan University, Wu Xi, China
\and
Suzhou University of Science and Technology, Su Zhou, China}
\maketitle  

%
\begin{abstract}
As global climate change and environmental issues escalate, carbon reduction has emerged as a paramount global concern. Agriculture accounts for approximately 30\% of global greenhouse gas emissions, making carbon reduction in this sector crucial for attaining global emission targets. Carbon sink trading serves as a supplementary mechanism to achieve carbon peaking and neutrality, helping to lower the rate of
carbon emissions. However, practical projects and research in the field of carbon sink trading are not enough currently. This work aims to thoroughly explore the cooperative models between farmers and retailers within the context of agricultural carbon sink trading, as well as the optimal decisions on the efforts to reduce carbon emission for both parties under different cooperative models. To this end, we delve into three distinct cooperative frameworks: the decentralized, the Stackelberg, and the centralized models, each accompanied by a corresponding differential
game model. The Hamilton-Jacobi-Bellman equation is utilized to investigate the equilibrium strategies of each participant under these three cooperative models, respectively. Furthermore, we conducte numerical simulations to analyze the carbon emission reduction efforts of farmers and retailers, the carbon emission reduction level of the agricultural supply chain, and the overall profits of the supply chain. We also compare scenarios with and without carbon sink trading to provide a comprehensive assessment. The numerical results indicate that the centralized model
excels in all aspects, followed by the Stackelberg model, with the decentralized model showing the weakest performance. Additionally, carbon sink trading can significantly increase the profits of the participants under each cooperative model. 

\end{abstract}
\newpage
\section{Introduction}
As the global economy continues to grow, the consumption of energy and resources, as well as emissions of greenhouse gases, exert substantial pressure on the ecological environment, leading to intensified climate change. The agricultural system, being one of the primary sources of greenhouse gas emissions, contributes more than 30\% to global emissions\footnote{https://www.chinanews.com.cn/gj/2021/03-11/9430282.shtml}. Consequently, reducing carbon emissions in agriculture is crucial for achieving global emission reduction goals. In this context, many countries have set targets for carbon emission reduction. 

Currently, countries worldwide primarily adopt carbon trading as a means to address carbon emissions issues, which typically occur in the midstream and downstream stages of the supply chain, such as processing and transportation. However, besides carbon trading, carbon sink trading is also important for the carbon market, serving as an effective measure for reducing agricultural carbon emissions with significant advantages. 
Agricultural carbon sink refers to the process of reducing greenhouse gas concentrations by absorbing carbon dioxide from the atmosphere through agricultural practices, such as crop planting and vegetation restoration. Enhanced carbon sequestration can slow the continued rise in atmospheric carbon dioxide concentrations and can aid in achieving carbon neutrality goals \cite{wan2021}. In carbon sink trading, farmers or agricultural enterprises can generate carbon credits or sinks through sustainable and low-carbon agricultural practices, which they can then sell to individuals or companies needing carbon sinks in the market. Carbon sink trading can effectively motivate farmers to reduce carbon emissions, thereby increasing carbon sequestration and reducing greenhouse gas emissions. Countries, such as the United States, Australia, and the European Union, have initiated agricultural carbon sink trading and achieved notable success.
In China, the nationwide carbon market is progressively being established. On January 22, 2024, the ceremony to launch the National Greenhouse Gas Voluntary Emission Reduction Trading Market was held in Beijing, marking the official restart of the china certified emission reduction (CCER) market since its suspension in 2017, presenting a key opportunity for the development of agricultural carbon sink projects in the country\footnote{http://paper.people.com.cn/zgnyb/html/2024-01/29/content\_26041024.htm}.

Simultaneously, increasing consumer awareness of low-carbon products has fueled their growing popularity in the market. Nowadays, consumers tend to prefer products with low-carbon properties when making purchasing decisions. The dual impetus of low-carbon policies and consumer demand preferences has incentivized carbon emission reduction in the agricultural product supply chain.  On the sales end, many large retailers have begun to set emission reduction targets and adopt  promotional strategies, such as Walmart's plan to reduce emissions by one billion tons, actively promoting a low-carbon supermarket image\footnote{https://finance.sina.com.cn/roll/2019-04-29/doc-ihvhiqax5710370.shtml}. On the production side, farmers can adopt emission reduction measures, such as reducing fertilizer usage and implementing conservation tillage, to increase soil carbon sequestration. They participate in carbon sink trading and earn additional income. The first carbon sink trading platform in Xiamen has been successfully established, facilitating the trade of 3,357 tons of carbon sink from agricultural tea gardens in two villages\footnote{http://www.fujian.gov.cn/xwdt/fjyw/202205/t20220506\_5903853.htm}.

However, the main obstacle to carbon reduction within the agricultural product supply chain is the imbalance in the distribution of revenues and costs related to carbon reduction among the participants. Participants including farmers and retailers, encounter different economic and technological constraints, as well as varying levels of awareness and willingness to invest in carbon reduction. So the obstacle lead to the challenge of promoting carbon reduction within the agricultural product supply chain. Specifically, the carbon emissions of agricultural products mainly occur in the upstream production phase. Farmers need to bear greater responsibility for carbon emission reduction and make high-cost investments. Still, these investments often cannot be converted into direct economic benefits in the short term, which undoubtedly increases the burden on farmers and reduces their enthusiasm for participating in carbon emission reduction. Retailers, as downstream participants in the supply chain, can directly benefit from consumer purchasing behavior and market demand, but they rely on upstream farmers to produce low-carbon agricultural products to meet market needs.

To address this issue and assist participants in making informed carbon reduction decisions within the agricultural product supply chain, this paper examines the impact of carbon sink trading and carbon emission reduction levels, constructing three differential game models: decentralized, Stackelberg, and centralized, for the supply chain system that includes the farmer and the retailer. Employing the Hamilton-Jacobi-Bellman equation, the optimal carbon reduction efforts and associated profits for each participant within these three modes of cooperation are explored. Our findings indicate that the centralized  model yields the best long-term benefits, with higher levels of carbon reduction and efforts compared to the decentralized and the Stackelberg models. In the Stackelberg game model, the retailer achieves higher overall profits by sharing costs with the farmer compared to the decentralized model. In addition, we found that carbon sink trading can increase the profit of the farmer and the retailer, incentivizing them to engage in carbon emission reduction.

\subsection{Related work}
Carbon trading is an effective means to promote carbon reduction. Some researches mainly focus on carbon quota allocation, carbon quota subsidies, and carbon emission reduction to explore the optimal carbon quota allocation mechanism and decision-making of supply chain members under carbon trading policies. 
Xu et al., \cite{xu2016}\cite{xu2017} investigated the production and emission reduction decision-making problem of 
Order-based supply chain, consisting of manufacturers and retailers under carbon cap and trade rules. 
Zhang et al., \cite{zhang2021} considered two scenarios, static carbon trading prices and dynamic carbon trading prices. They developed an evolutionary game model involving the government and manufacturers to investigate how government policies influence manufacturer decisions and the carbon trading market. 
Cai et al., \cite{cai2023} studied a supply chain model consisting of a supplier and a manufacturer, and discussed the optimal pricing and carbon reduction decisions of the supply chain members. In the context of carbon sink trading, Wang et al., \cite{wang F2021} considered carbon sinks can fully utilize the characteristics of natural ecosystems and will be an important means in the process of carbon emission reduction. 
Ke et al., \cite{ke2023china} conducted research on forest carbon sink trading. The study assessed the potential of forest carbon sink trading in China. Wang et al., \cite{wang2015} addressed the gap in carbon sink trading in China and, based on a duopoly model, examined the impact of carbon sink mechanisms on the profits of emission trading participants and industry output. 
The results indicate that when companies exceed their carbon emission limits, they can purchase carbon sinks to compensate for the excess emissions. 

In terms of supply chain coordination, many researches have started applying supply chain contracts to agricultural product supply chains.
Yu et al., \cite{yu2017} proposed a Stacklberg game model to coordinate pricing and service levels in a three-stage Agricultural Product Supply Chain that includes third-party logistics service providers. 
The model takes into account both supplier-led and logistics service provider-led scenarios. 
Furthermore, Yan et al., \cite{yan2017}examined the optimal decision-making challenges in an Internet of Things (IoT) fresh agricultural product supply chain, which comprises manufacturers, distributors, and retailers. 
They proposed improved revenue-sharing contracts to achieve supply chain coordination. Additionally, Song et al., \cite{song2019} investigated the decision-making problems in a three-stage supply chain, including fresh e-commerce companies, third-party logistics service providers, and community convenience stores, under both centralized and decentralized model. 
They introduced cost-sharing contracts and profit-sharing contracts to achieve coordination and maximize the profits of all parties in the supply chain. 
Ma et al., \cite{ma2019} studied the coordination issues in a three-stage supply chain system consisting of a supplier, a third-party logistics service provider, and a retailer supplying seasonal fresh agricultural products to customers. 

In practice, carbon emissions reduction in the agricultural product supply chain is a long-term process \cite{zu2018}, where the outcomes of previous stage decisions can impact the subsequent carbon reduction decisions, and the carbon emissions of the final product are influenced by the carbon reduction decisions made by different participants in the supply chain. 
Therefore, it is necessary to analyze it from a dynamic perspective. 
Differential game theory is the useful tool to study the dynamic decision problems involving continuous time and continuous actions, in which participants aim to achieve optimal decisions through dynamic optimization. 
Some researches have employed the approach of differential game theory to investigate carbon emissions reduction in supply chains.
Wang et al., \cite{wang2021} developed three differential game models for a supply chain system consisting of two suppliers and a single manufacturer, examining the impact of different carbon emission allocation rules on carbon reduction in the supply chain. 

Some studies have explored supply chain decision-making under carbon trading \cite{xu2016,xu2017,zhang2021,cai2023}. 
Others have examined the significance of carbon sink trading, primarily focusing on carbon sink assessment and the factors influencing carbon sink trading \cite{ke2023china,wang2015}. 
Notably, several of these studies have been applied in the forestry sector.
However, research on agricultural carbon sink trading, particularly in terms of coordinated decision-making within the supply chain, is 
 not too enough. Also, previous research on low-carbon supply chains has primarily focused on carbon emissions during processing and transportation. 
To enrich the research of carbon sink trading, we mainly investigate carbon emissions during the production of agricultural products and innovatively incorporate carbon sink trading into the agricultural supply chain. We utilize the differential game approach to examine the carbon emission reduction decisions of farmers and retailers within the supply chain.

\subsection{Organization}
The remaining part of this paper is organized as follows: Section 2 provides a detailed description of the model framework along with the necessary assumptions. 
In Section 3, we introduce and explain three different cooperation models, namely decentralized decision model, Stackelberg game model, and centralized decision model, and analyze their respective equilibrium solutions. 
Numerical analysis and sensitivity analysis are conducted in Section 4. 
Finally, Section 5 presents the conclusion, discussion, and outlines potential directions for future research.

\section{Problem Description and and Assumptions}

\subsection{ Description of the problem}

Due to the widespread distribution of farmers, it is difficult to account for the carbon sink produced by an individual farmer.
Therefore, in practice, it is common to aggregate the overall carbon sinks produced by agricultural activities at the village or township level for the purposes of unified accounting and trading.
To simplify the expression and calculation, we treat the farmers of the entire village as a single collective entity, denoted by $f$.
This work focuses on a low-carbon agricultural product supply chain, consisting of one farmer $f$ and one retailer $r$. 
The primary responsibility of the farmer is to undertake efforts such as implementing advanced production techniques and enhancing soil quality to achieve the carbon emission reduction targets.
In detail, they adopt production technologies such as nitrogen oxide emission reduction in farmland to enhance production efficiency and reduce the emission intensity per unit of product. 
Additionally, they improve soil quality to enhance carbon sequestration capacity in farmland, including implementing conservation tillage, promoting straw returning to the field. 
On the other hand, the retailer enhance the environmental image of agricultural product brands through active low-carbon promotion and advertising activities, thereby attracting more consumers. 
These combined efforts aim to enhance the brand's environmental image, collectively referred to as low-carbon promotional initiatives.
Figure 1 demonstrates the process of carbon emission reduction in the agricultural product supply chain.
\setcounter{figure}{0} 
\begin{figure}[H]
  \centering
 \includegraphics[width=0.95\textwidth]{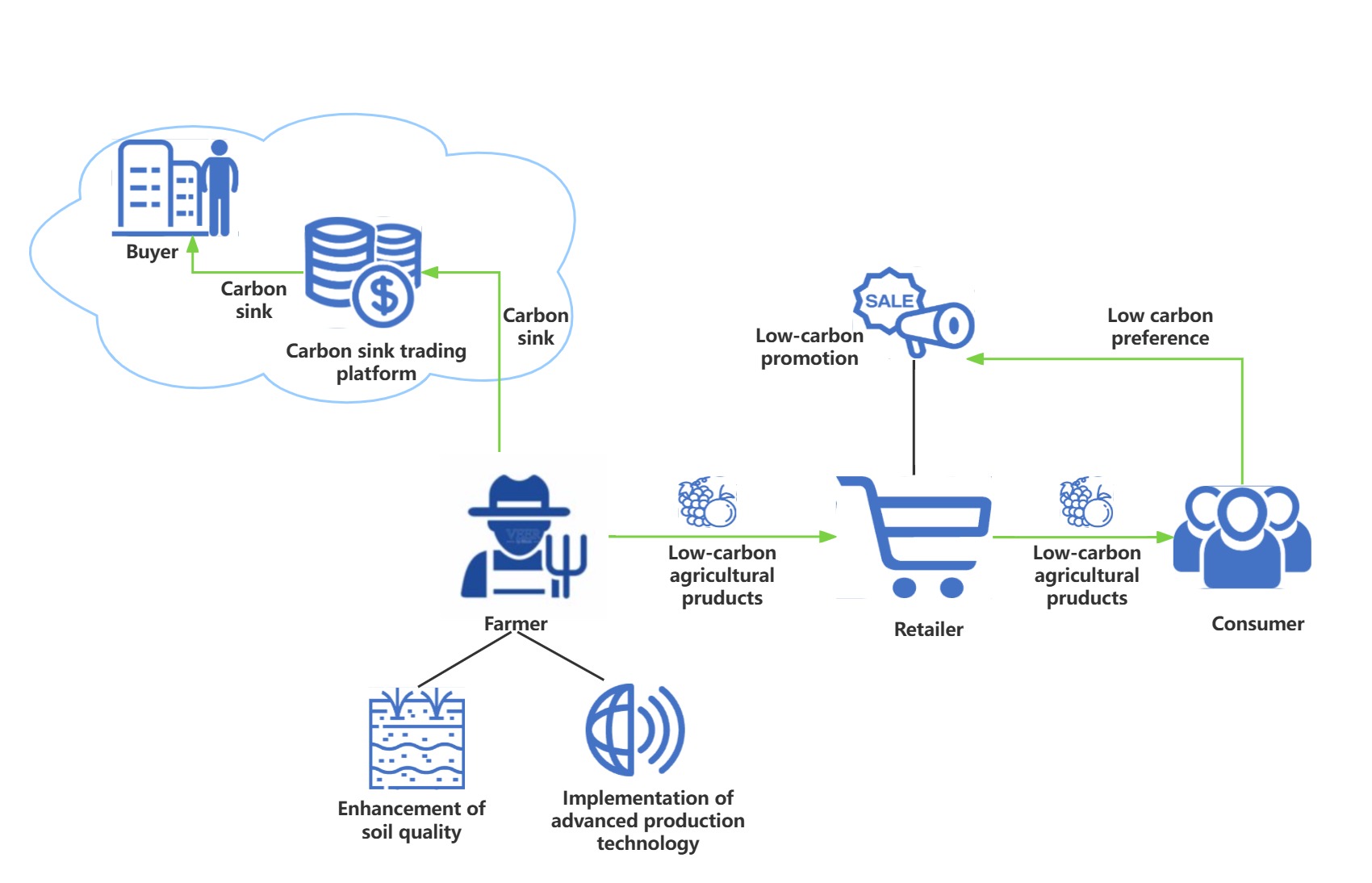}
  \caption{The process of carbon emission reduction in the agricultural product supply chain}
   \label{process}  
\end{figure}

\subsection{Problem assumptions}

The following assumptions are necessary for our discussion.

\noindent {\bf Assumption 1. }
Generally, products that consume less energy and emit lower amounts of carbon dioxide and other greenhouse gases during their production and consumption processes are referred to as low-carbon products. 
Given consumers' preference for low-carbon products, the level of carbon emissions reduction $H$ in the agricultural supply chain has become a critical factor that influences the marketing of agricultural products and consumer purchasing decisions. This level $H$ is influenced by two factors: the emission reduction and carbon sequestration efforts $E_f(t)$ of the farmer and low-carbon promotion efforts $E_r(t)$ of the retailer. Similar to \cite{wang2021}, we can describe the dynamic process of $H(t)$ in agricultural products using the following differential equation:
\small{
\begin{equation}
\begin{aligned}
\dot{H}(t)=\mu_f E_f(t)+\mu_r E_r(t)-\delta H(t), 
\end{aligned}
\end{equation}}
where $\mu_f$ and $\mu_r$ represent the impact coefficient of the farmer's and the retailer's efforts to reduce carbon emission, on the level of carbon emission reduction, and $\delta>0$ represents the decay coefficient of the carbon emission reduction level. 
 
\noindent {\bf Assumption 2. }
Faced with environmental deterioration, people's awareness of low-carbon living has risen, leading more consumers to pay more for low-carbon products \cite{wu2020}. Additionally, as the carbon emission reduction level of the supply chain improves, the low-carbon attributes of products become increasingly significant, leading to a larger market demand.
Also, different cooperation models may lead to information asymmetry between  farmers and retailers.
This asymmetry may cause the supply of agricultural products to not match market demand.
As a result, retailers could face stock shortages.
Alternatively, they might experience issues with inventory overstock. 
Therefore, we assume that the supply $Q(t)$ and demand $D(t)$ of agricultural products are influenced by the level of carbon emission reduction $H(t)$, the price per unit of agricultural product $p$, and consumers' preferences for low-carbon products $\theta$.
Following \cite{cai2023}, 
the supply and demand function of agricultural products are 
\small{\begin{eqnarray}
Q(t)&=&(Q_0+ap)\theta H(t), \\
D(t)&=&(D_0-bp)\theta H(t), 
\end{eqnarray}}
where $Q_0\geq0$ and $D_0\geq0$ represent the initial supply and initial demand of agricultural products, respectively; $a\geq0$ and $b\geq0$ represent the influence of agricultural product prices on product supply and demand, respectively. 

\noindent {\bf Assumption 3. }
The carbon emission reduction during the agricultural production process is crucial for the carbon emission reduction of the agricultural product supply chain. 
By implementing advanced technologies and enhancing soil quality, farmers can lower carbon emissions, and thus increase carbon sink $F$, a distinctive feature of the agricultural industry.
Assuming that the carbon emission reduction level of agricultural products is correlated with the farmer's efforts $E_f(t)$ in emission reduction and carbon sequestration, represented by the coefficient $\omega$ \cite{wang2021}. Similar to \cite{xia2017,po2007}, we assume that the farmer's initial carbon sink per unit product is 1. Thus, the total carbon sink $F(t)$ at time $t$ can be expressed as:
\small{\begin{eqnarray}
   F(t)=(1+\omega E_f(t))Q(t), 
\end{eqnarray}}
where $Q(t)$ is the supply of the farmer at time $t$. 

\noindent {\bf Assumption 4. } The costs of emission reduction efforts by the farmer and the retailer have a quadratic relationship with $E_f(t)$and $E_r(t)$, respectively. The form of quadratic function, which is widely adopted to describe the cost pattern in the literature
\cite{Ji2017,Yang2023}. To simplify the problem and highlight the research focus, we assume that the production cost of agricultural products is zero\cite{wang2021}. Therefore, the costs ofcarbon emission reduction efforts by the farmer and the retailer are
\small{\begin{equation}
C_s(t)=\frac{1}{2} \lambda_f E_f^2(t),~~ 
C_m(t)=\frac{1}{2} \lambda_r E_r^2(t),  
\end{equation}}
where $\lambda_f$, $\lambda_r$ represent the effort cost coefficients for the farmer and the retailer to reduce carbon emissions, respectively. 

Table~\ref{tab1} provides the notations used in this paper.

\small{\FloatBarrier
\begin{table}[!ht]
\addtocounter{table}{0}
 \vspace{-5pt}
 \renewcommand{\tablename}{Table}
\captionsetup{justification=raggedright, singlelinecheck=false} 
  \caption{Major notations}
  \label{tab1}
   \begin{tabularx}{340pt}{cX}
\toprule
             \textbf{Notations}                   &    \textbf{Explanation}
\\ \hline
     $t$       &  Time period
\\
       $H(t), H(0)$   & Carbon emission reduction level  of agricultural products at time $t$, and initial carbon emission reduction level, $H(0) \geq 0$. 
\\
        $E_f(t), E_r(t)$    &  The  emission reduction and carbon sequestration efforts of the farmer at time $t$, the low-carbon promotion efforts $E_r(t)$ of the retailer at time $t$.  
\\
        $\lambda_f, \lambda_r$          & Effort cost coefficients of  the farmer and the retailer to reduce carbon emissions, $\lambda_f, \lambda_r> 0$. 
\\
      $\mu_f, \mu_r$       & The impact coefficient of the farmer's and the retailer's carbon emission reduction efforts on the level of carbon emission reduction, $\mu_f, \mu_r>0$. 
\\
     $\omega$   & The impact coefficient of the farmer's efforts to reduce carbon emissions on the carbon emission level, $\omega>0$. 
\\
     $p_f, p_r, p, p_c$  & Marginal profits of the farmer and the retailer, unit agricultural product price, unit carbon sink price, $p_f, p_r, p, p_c\geq0$. 
\\
    $Q(t), D(t)$    & Supply function for the farmer at time $t$, and demand function for the retailer at time $t$. 
\\
     $\theta$  & Consumers’ low-carbon preferences, $\theta>0$. 
\\
      $ a, b $   & The influence of agricultural product prices on product supply and demand, $ a>0, ~b>0$. 
\\
      $ F(t)$  &The total carbon sink at time $t$. 
\\
       $ \delta$  & Attenuation coefficient of carbon emission reduction level, $ \delta>0$.  
\\
     $\rho$ & profit discount rate, 
     $\rho>0$. 
\\
\bottomrule
\end{tabularx}
\end{table}}

\section{Analysis for Differential Game Equilibrium}

Building on the assumptions mentioned earlier, we have developed three differential game models: decentralized, Stackelberg, and centralized. These game models are used to analyze the equilibrium strategies of the farmer and the retailer across three different cooperative scenarios in this section. For simplicity, the parameter $t$ will be omitted in the following sections.

\subsection{Decentralized decision-making model}
Under the decentralized decision-making model, indicated by the superscript GD, both the farmer $f$ and the retailer $r$ aim to maximize their respective profits over an infinite time horizon, with a profit discount rate assumed to be $\rho$. Thus, the objective functions for the farmer and the retailer are formulated as (6) and (7). 
\small{\begin{align}
\max J_{f}^{GD} & =\int_0^{\infty} e^{-\rho t}\left[p_f Q  + p_c (1+\omega E_f)Q -\frac{1}{2} \lambda_f E_f^2\right]dt, \label{qf1}\\
\max J_r^{GD} & =\int_0^{\infty} e^{-\rho t} \left[p_r D -\frac{1}{2} \lambda_r E_r^2 \right]dt.  \label{qr1} 
\end{align} }

Proposition \ref{Proposition1} propose the equilibrium strategies for both parties involved in the decentralized decision-making model. The complete proof of this proposition can be found in the full version, due to the limit of the space.
\begin{Proposition}\label{Proposition1}
Under the decentralized decision-making model, the equilibrium carbon emission reduction efforts of the farmer and the retailer are
\small{\begin{equation*}
    \begin{aligned}
 {E_f^{GD}}^* =\frac{p_c \theta\omega k_1 H + \mu_f(2A^{GD} H +B^{GD})}{\lambda_f}, ~~
{E_r^{GD}}^* =\frac{\mu_r M^{GD}}{\lambda_r}, 
\end{aligned}   
\end{equation*}}
the equilibrium level of carbon emission reduction is
\small{ \begin{equation*}
 \begin{aligned}
H_d^{GD} &=\frac{\lambda_r\mu_f B^{GD}+\lambda_f\mu_r M^{GD}}{\lambda_f\lambda_r \delta-2A^{GD}\mu_s\lambda_r-p_c\omega k_1 \mu_f\lambda_r},
\end{aligned}   
\end{equation*}}
the optimal profits of the farmer and the retailer are
\small{\begin{equation*}
    \begin{aligned}
{V_f^{GD}}^*= A^{GD} H^2 + B^{GD} H + C^{GD}, ~~
{V_r^{GD}}^*= M^{GD} H + N^{GD},    
  \end{aligned}
\end{equation*}}
where $A^{GD}, B^{GD}, C^{GD}, M^{GD}, N^{GD}, T^{GD}$are shown in \eqref{q1}
\small{\begin{equation}
  \begin{aligned}  
A^{GD}&= \frac{2\rho \lambda_f+4\lambda_f \delta-4\mu_f p_c \omega k_1- \sqrt{ \triangle^{GD} }}{8 \mu_f^2}, \\
M^{GD}&=\frac{p_{m} k_{2} \lambda_{s}}{p_{c} \omega k_{1}\mu_{s}+2A^{GD}\mu_{s}^{2}+\lambda_{s}\rho-\delta\lambda_{s}}, \\
B^{GD}&=-\frac{\left(2 \lambda_r^{2} A^{GD}M^{GD}+k_1 \lambda_r \mathit{p_c}+k_1\lambda_r\mathit{p_f}\right) \lambda_f}{\lambda_r \left(\mathit{p_c} \omega k_1 \mathit{\mu_f}+2 A \, \mathit{\mu_f}^{2}-\lambda_f \rho-\delta \lambda_f\right)}, \\
C^{GD} &=\frac{B^{GD} \left(\lambda_{s}^{2}\lambda_{m} B^{GD}+2M^{GD}\lambda_{m}^{2} \lambda_{s}\right)}{2 \lambda_{m}\lambda_{s}\rho}, \\
N^{GD}&=\frac{M^{GD} \left(2 \lambda_{m}\mu_{s}^{2} B^{GD} +M^{GD}\mu_{m}^{2}\lambda_{s}\right)}{2 \lambda_{m}\lambda_{s}\rho},
  \end{aligned}
  \label{q1}
\end{equation}}
with $\triangle^{GD} = (4 \mu_f  p_c \omega k_1-4\lambda_f \delta-2\rho\lambda_f)^2- (4\mu_f p_c\omega k_1 )^2 > 0$, $k_1=(Q_0+ap)\theta$, and $k_2=(D_0-bp)\theta$. 
\end{Proposition}

\subsection{Stackelberg game model}
Under this model, the retailer provides subsidies to the farmer to encourage efforts towards reducing carbon emissions.
In the first stage, the retailer acts as a leader, determining their efforts to low-carbon promotion $E_r$ and the subsidy ratio $x_f$ for the farmer. Then, in the second stage, acting as a follower, the farmer makes her optimal responses based on the decisions made by the retailer in the first stage. Thus, the objective functions for farmer and the retailer are formulated as follows:
    \begin{align}
\max J_r^{GS}&=\int_0^{\infty} e^{-\rho t} \bigg[p_r D-\frac{1}{2} \lambda_r E_r^2-\frac{1}{2}x_f\lambda_f E_f^2 \bigg]dt,\label{qf2}\\
\max J_{f}^{GS} & =\int_0^{\infty} e^{-\rho t}\left[p_fQ+p_c (1+\omega E_f)Q -(1-x_f)\frac{1}{2} \lambda_f E_f^2\right]dt.\label{qr2}
\end{align}
\label{jsm}

Similarly, the equilibrium strategies for both parties under the Stackelberg game model are provided in Proposition \ref{Proposition2}, whose proof is also in the full version.
\begin{Proposition}\label{Proposition2}
Under the Stackelberg game model, the equilibrium carbon emission reduction efforts of the farmer and the retailer are 
   \small{
    \begin{equation*}
       \begin{aligned}
 {E_f^{GS}}^* =\frac{p_c \omega k_1 H + \mu_f(2A^{GS} H+B^{GS})}{\lambda_f}, ~~{E_r^{GS}}^* =\frac{\mu_r (2M^{GS}H+N^{GS})}{\lambda_r}, \\  
\end{aligned}   
\end{equation*}}
the equilibrium level of carbon emission reduction is
  \small{
  \begin{equation*}
       \begin{aligned} 
H_d^{GS}&=\frac{2(\lambda_f\rho\mu_r +\lambda_r\mu_f)N^{GS}+\lambda_r\mu_f B^{GS}}{4\lambda_r\mu_f M^{GS}-\rho \lambda_r p_c \omega k_1-4\lambda_f\mu_r\rho M^{GS}-2\lambda_r\mu_f A^{GS}+2\lambda_r\lambda_f\rho\delta},
\end{aligned}   
\end{equation*}}
the equilibrium  subsidy ratio  is
  \small{
  \begin{equation*}
       \begin{aligned} 
x_f^{{GS}^*}&=\frac{(4M-2A)H+2N-B-p_c\omega k_1 H/\mu_f}{(4M+2A)H+2N+B+p_c\omega k_1 H/\mu_f},
\end{aligned}   
\end{equation*}}
The optimal profits for  the farmer and the retailer are
\small{
\begin{equation*}
\begin{aligned}
{V_f^{GS}}^*= A^{GS} H^2 + B^{GS} H + C^{GS}, ~~{V_r^{GS}}^*= M^{GS}H^2+N^{GS} H + F^{GS},    
\end{aligned}
\end{equation*}}
where $A^{GS}, B^{GS}, C^{GS}, M^{GS}, N^{GS}, T^{GS}$are shown in(\ref{q2})
\small{
\begin{equation}
  \begin{aligned}  
A^{GS}&=\frac{2\lambda_r\lambda_f\rho \delta-p_c\omega\ k_1\lambda_r\rho\mu_f-\lambda_r\lambda_f\rho^3-2\mu_f^2\lambda_r M^{GS}-4\mu_r^2\lambda_f M^{GS}-\sqrt{\triangle^{GS1}}}{2\mu_f^2\lambda_r}, \\
M^{GS}&=\frac{2\lambda_f\lambda_r\rho\delta-p_c\omega k_1\rho\lambda_r\mu_f-\lambda_f\lambda_r\rho^3-2A^{GS}\lambda_r\mu_f^2-\sqrt{\triangle^{GS2}}}{4(\lambda_r\mu_f^2+\lambda_f\mu_r^2)}, \\
B^{GS}&=\frac{\lambda_f\lambda_r((p_c +p_f)k_1+p_r k_2)}{\lambda_f\lambda_r(\rho-\delta)-p_c\omega\lambda_r\mu_f k_1-2A^{GC}(\lambda_r\mu_f^2+\lambda_f\mu_r^2)}, \\
C^{GS}&=\frac{\mu_f^2\lambda_r{B^{GS}}^2+(4\mu_r^2\lambda_f+2\mu_f^2\lambda_r)B^{GS}N^{GS}}{4\rho^3\lambda_f\lambda_r}, \\
N^{GS}&=\frac{2\lambda_r\mu_f^2B^{GS}(A^{GS}+2M^{GS}-4\lambda_r\lambda_f p_r\rho^2 k_2-p_c\omega k_1\lambda_r\mu_f\rho B^{GS})}{4\rho^3\lambda_f\lambda_r-2p_c\omega k_1\mu_f\rho\lambda_r+4\lambda_f\lambda_r\delta\rho+4\mu_f^2\lambda_r(A^{GS}+2M^{GS})+8\mu_r^2\lambda_f M^{GS}}, \\
F^{GS}&=\frac{\mu_f^2\lambda_r({B^{GS}}^2+ 4N^{GS}B^{GS})+4(\lambda_r\mu_f^2+\lambda_f\mu_r^2){N^{GS}}^2}{8\rho^3\lambda_r
\lambda_f},
  \end{aligned}
  \label{q2}
\end{equation}}
with $k_1=(Q_0+ap)\theta$, $k_2=(D_0-bp)\theta$,\\
\begin{equation}
    \begin{aligned}
\triangle^{GS1}&=2(\mu_f \omega k_1 p_c\rho)^{2} +2k_1 p_c\omega\lambda_f\mu_f \rho^{2}(\rho^{2}-2\delta) +\lambda_f\rho^{3}(\lambda_f\rho^{3}+4M\mu_f^{2}-4\lambda_f\rho\delta)\\
&+4M^{2}\mu_f^{4}-8M\lambda_f\rho\delta\mu_f^{2}+4\lambda_f^{2}+4(\lambda_f\rho\delta)^{2}+8M\mu_r^{2}\lambda_f(\lambda_r\rho\omega p_c k_1+\lambda_f\lambda_r\rho^{3}\\
&+2M\mu_f^{2}\lambda_r+2M\mu_r^{2}\lambda_f-2\lambda_r\lambda_f\rho\delta)>0, and\\
\triangle^{GS2}&=2\lambda_r(k_2 p_c\rho\omega)^{2}(\lambda_r\mu_f^{2}+\lambda_f)+2k_2\lambda_r\lambda_f p_c \rho\omega(\lambda_r\rho^{3}\mu_f-2A\mu_r^{2}-2\lambda_r\rho\delta)\\
&+4\lambda_r\lambda_f\rho^{2}(\rho^{3}+\rho^{2}\delta+\delta^{2})+4A\lambda_r\lambda_f\mu_f^{2}(\lambda_r\rho^{3}-A\mu_r^{2}-2\lambda_r\rho\delta)>0.\nonumber
    \end{aligned}
\end{equation}
\end{Proposition}

\subsection{Centralized decision-making-making model}
Under this model, the cooperation among the farmer and the retailer are aimed at maximizing the total profits of the entire system. 
Unlike the previous two models, the objective function in this case is derived by discounting the combined profits of the two participants, with a profit discount rate assumed to be $\rho$.
Thus, we propose the objective function in the centralized decision-making model.
    \begin{align}
\max J_{fr}^{GC}& =\int_0^{\infty} e^{-\rho t}\bigg[p_f Q+p_r D + p_c (1+\omega E_f)Q-\frac{1}{2}\lambda_f E_f^2 -\frac{1}{2}\lambda_r E_r^2\bigg] dt. \label{qfr}
 \end{align}

Under the centralized decision-making model, the equilibrium strategies of the corresponding differential game are proposed in the following Proposition. We also move its proof to the full version.
\begin{Proposition}\label{Proposition3}
Under the centralized decision-making model, the equilibrium carbon emission reduction efforts of the farmer and the retailer are 
\begin{equation*}
    \begin{aligned}
{E_f^{GC}}^* &=\frac{\mu_f(2A^{GC}H +B^{GC})}{\lambda s}, ~~ 
{E_r^{GC}}^* =\frac{p_c \theta \omega \theta(D_0-bp) H + \mu_r(2A^{GC} H+B^{GC} )}{\lambda_r}, \\ 
    \end{aligned}
\end{equation*}
the equilibrium level of carbon emission reduction is
\small{
\begin{equation*}
    \begin{aligned}
H_d^{GC}&=\frac{(\lambda_f\mu_r+\lambda_r\mu_f)B^{GC}}{\lambda_f\lambda_r\delta-\lambda_r p_c\omega k_1-2(\lambda_f\mu_r+\lambda_r\mu_f)A^{GC}},
    \end{aligned}
\end{equation*}}
 the optimal profit of the supply chain is
\small{\begin{equation*}
    \begin{aligned}
  V_{fr}^{{GC}^*}&= A^{GC} H^2 + B^{GC} H + C^{GC},      
    \end{aligned}
\end{equation*}}
where $A^{GC}, B^{GC}, C^{GC}$are shown in (\ref{q3})
\end{Proposition}
\small{
\begin{align}       
A^{GC}&=\frac{2\delta\lambda_f-\rho\lambda_f-2\lambda_r\mu_fp_c\omega k_1-\sqrt{\triangle^{GC}}}{4\lambda_r \mu_f^2+4\lambda_f \mu_r^2}, \nonumber\\
B^{GC}&=\frac{\lambda_f\lambda_r((p_c +p_f)k_1+p_r k_2)}{\lambda_f\lambda_r(\rho-\delta)-p_c\omega\lambda_r\mu_f k_1-2A^{GC}(\lambda_r\mu_f^2+\lambda_f\mu_r^2)}, \nonumber\\
C^{GC}&=\frac{B^{{GC}^{2}}(\lambda_r\mu_f^2+\lambda_f\mu_r^2)}{2\lambda_r\lambda_f\rho},\label{q3}
\end{align}}
with $\triangle^{GC} =4\lambda_r(p_c \omega k_1)^2(2\lambda_r\mu_f^2+\lambda_f\mu_r^2+\lambda_r^2\lambda_f(\rho-2\delta)(4p_c\omega k_1\mu_f+\rho-2\delta))$, \\
$k_1=(Q_0+ap)\theta$, 
and $k_2=(D_0-bp)\theta$. \\

\section{Numerical simulation and analysis}
Referring to \cite{zhang2021,cai2023},  the relevant parameter in our numerical experiments are set as follows:
$\lambda_f = 500, 
\lambda_r= 200, 
\mu_f = 1. 5, 
\mu_r = 0. 5, 
\omega= 0. 4, p_f = 5, 
p_r = 10, 
p=25, 
p_c = 0. 5, 
a = 3, 
b = 2, 
\delta= 1, 
\rho = 0. 7, 
\theta = 0. 8, 
Q_0 = 300, 
D_0 = 250. $
The outcomes of the decision model are illustrated using distinct colors to distinguish among different cooperation modes. Specifically, red represents the decentralized decision-making model, green denotes the Stackelberg game model, and blue indicates the centralized game model. It's important to note that solid lines symbolize scenarios with carbon sink trading, whereas dashed lines depict scenarios without carbon sink trading.

\subsection{Impact of time $t$ and parameter $\omega$ on the level of carbon emission reduction and total profit}
This section discusses the changes in the level of carbon emission reduction and total supply chain profits over time, considering three types of cooperation models. 

\begin{figure}[H]
  \centering
  \subfigure[The level of carbon emission reduction under three modes]{\includegraphics[width=0.49\textwidth]{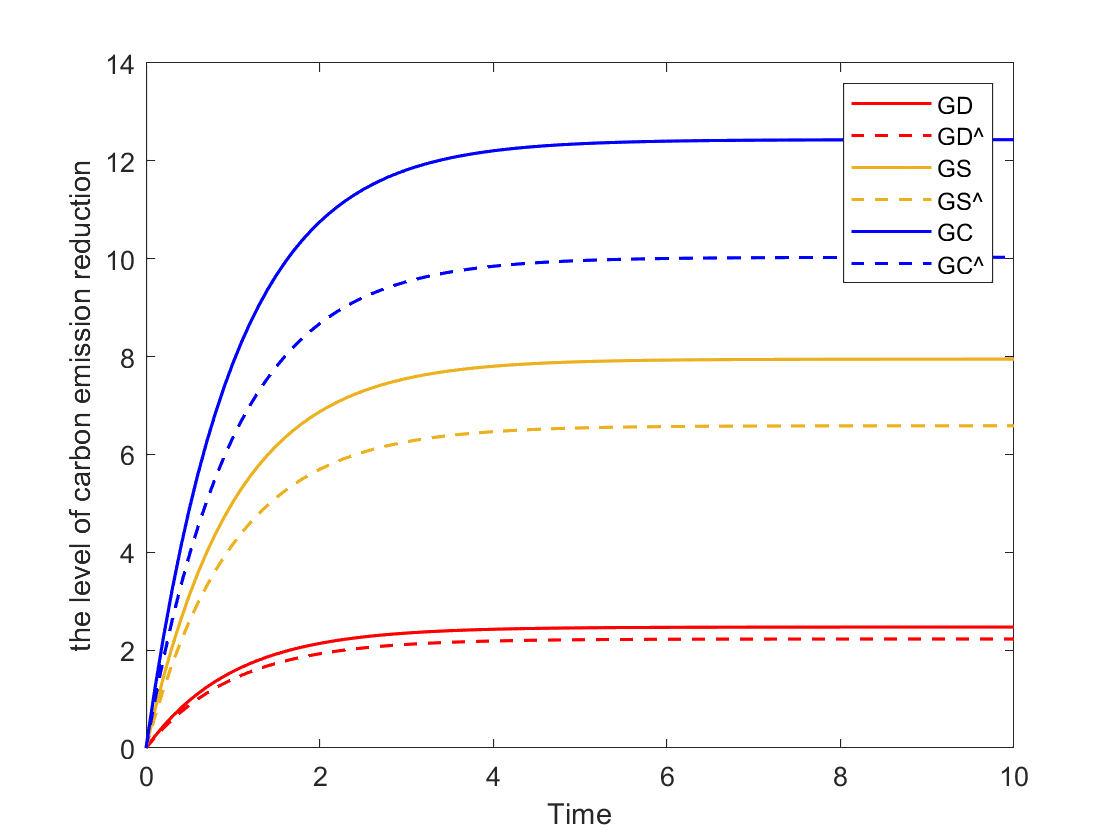}}
  \subfigure[Total profit under three modes]{\includegraphics[width=0.49\textwidth]{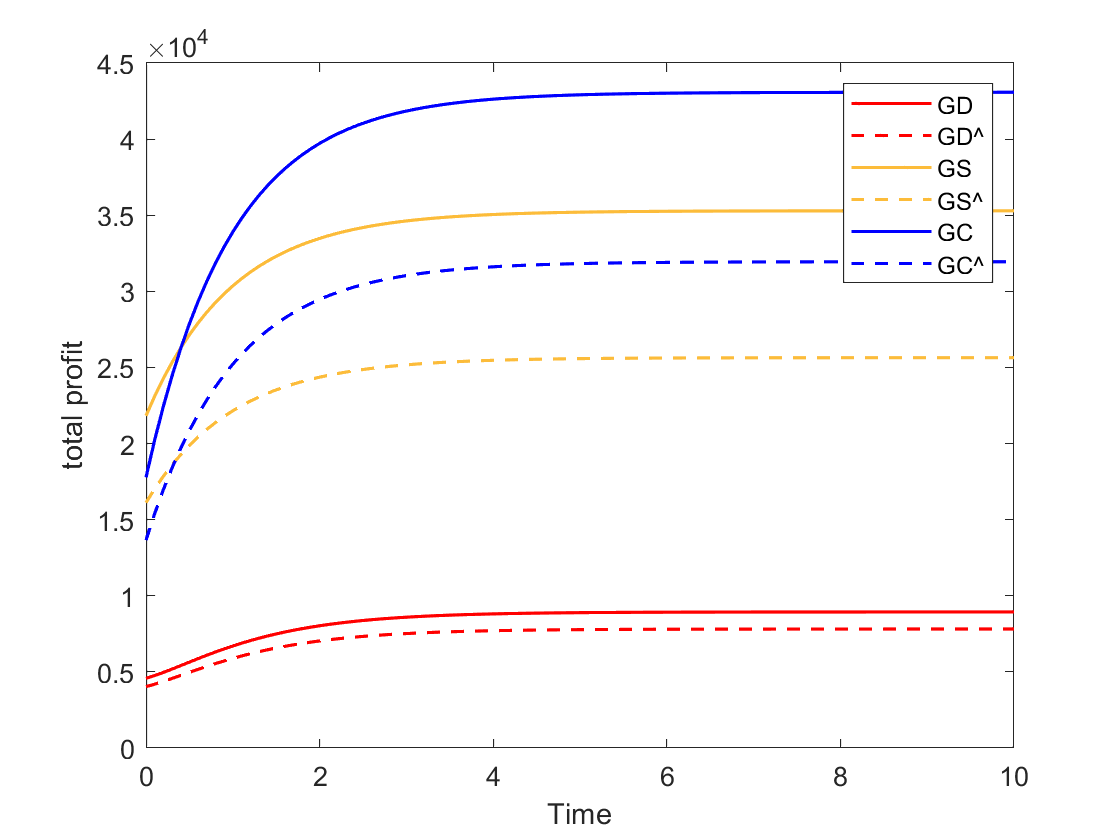}}
    \caption{The level of carbon emission reduction and total supply chain profits}
    \label{carbon emission}  
\end{figure}

Fig.\ref{carbon emission}-a depicts the carbon emission reduction level trajectories under three types of cooperation models, which increase over time $t$ and stabilize as $t$ approaches infinity. Clearly, the centralized decision-making model achieves the highest level of carbon emission reduction among the three models, followed by the Stackelberg game model, while the decentralized decision-making model exhibits the lowest level of carbon emission reduction. 

Fig.\ref{carbon emission}-b depicts the changes of the total profit in the supply chain over time, under three types of cooperation models. 
We can clearly observe that the decentralized decision-making model shows the lowest total profit in the supply chain. While the Stackelberg game model initially outperforms the centralized decision model in terms of profit, over time, the centralized model gradually reveals more substantial advantages and significantly exceeds the other two models.
 
Specifically, the centralized model consolidates carbon emission reduction resources from stakeholders like the farmer and the retailer to optimize resource allocation. In contrast, the decentralized model, lacking unified coordination, risks redundant investments and resource wastage. In the Stackelberg model, the retailer set emission standards and targets, guiding the farmer and driving the supply chain towards a low-carbon transformation.
Additionally, scenarios involving carbon sink trading outperform those without carbon sink trading under three types of cooperation models, indicating that appropriate incentive mechanisms can impove both the carbon reduction efficiency and the profits of supply chain participants.

\subsection{Impact of time $t$ and parameter $\omega$ on the farmer's and the retailer's profits}
This section discusses the changes in the profits of the farmer and the retailer over time, under the decentralized decision model and the Stackelberg game model.

\begin{figure}[H]
  \centering
  \subfigure[Farmer's profit under different modes]{\includegraphics[width=0.49\textwidth]{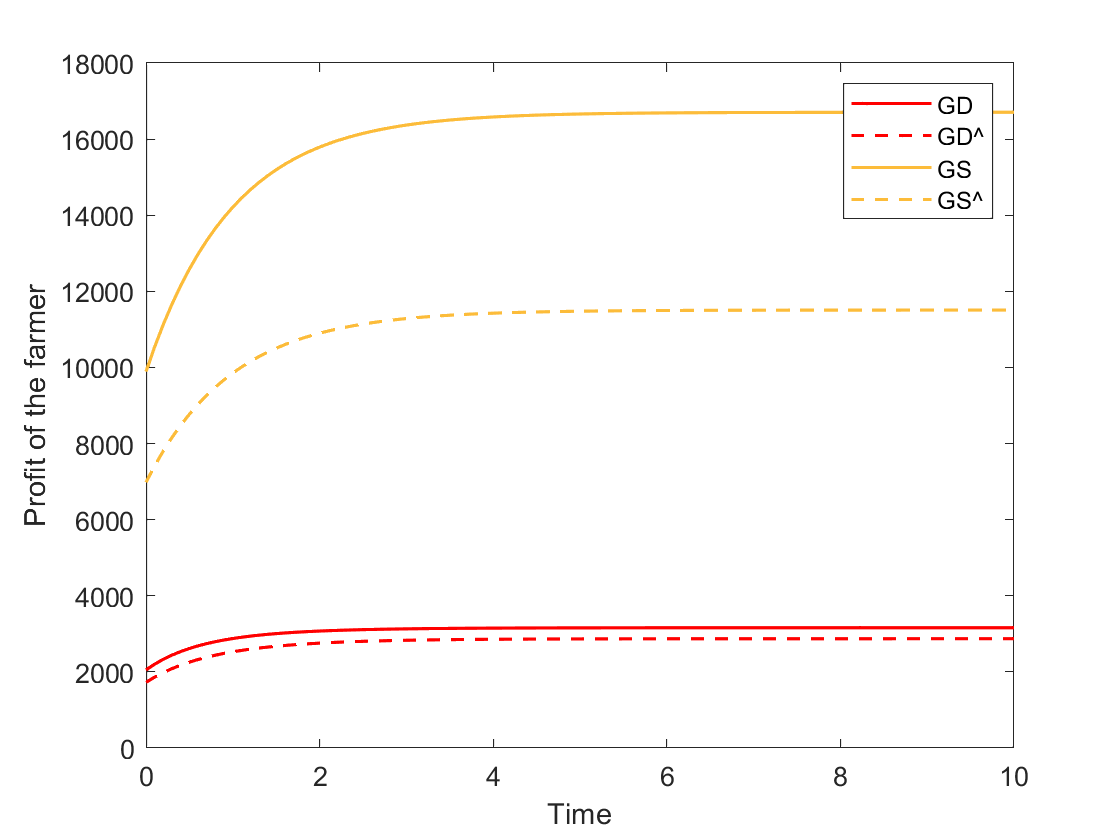}}
  \subfigure[Retailer's profit under different modes]{\includegraphics[width=0.49\textwidth]{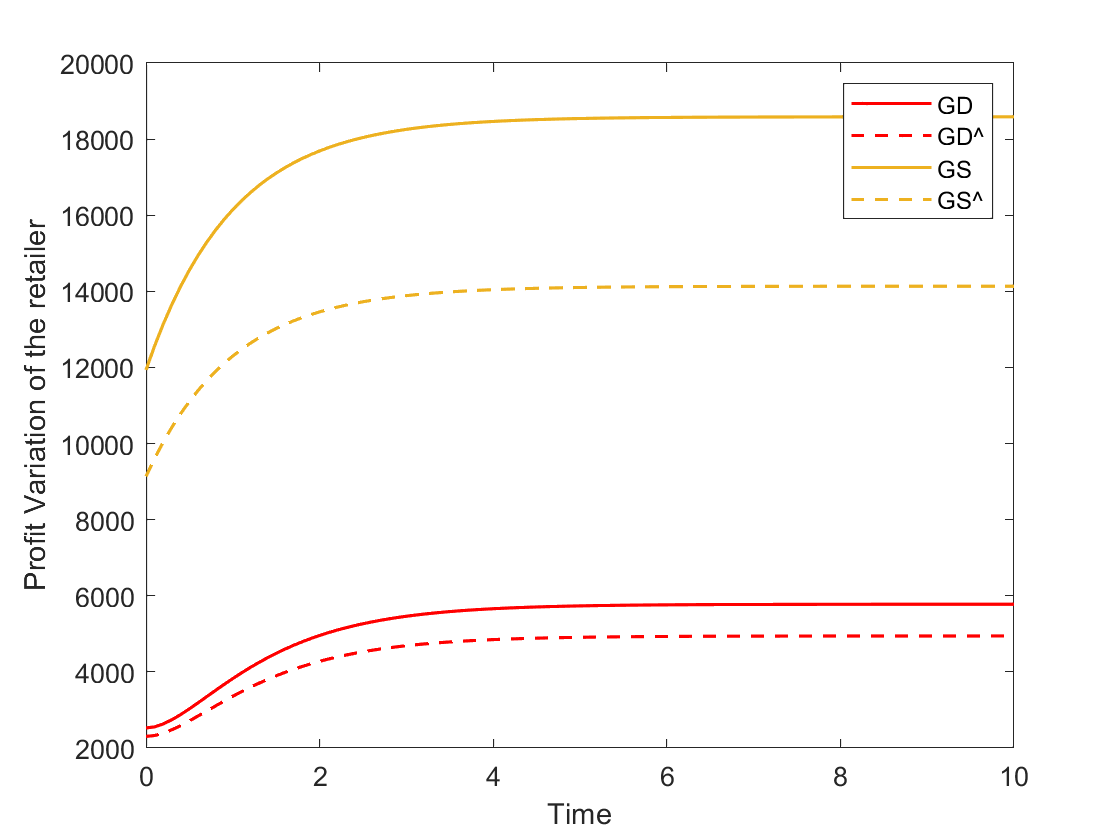}}
  \caption{Profits of farmer and retailer under  decentralized and Stackelberg  model.}
  \label{profit}
\end{figure}

Fig.\ref{profit}-a depicts the changes of the farmer's profits over time, under two models.
And Fig.\ref{profit}-b depicts the changes of the retailer's profits over time, under two models. From the two figures, it is evident that under the Stackelberg model, both the farmer and the retailer achieve higher profits compared to those in the decentralized model. 

In the Stackelberg model, the retailer provides partial cost subsidies to the farmer, thereby encouraging them to adopt carbon reduction measures. On the other hand, in the decentralized model, the farmer bears greater costs for carbon reduction. Additionally, the absence of effective communication between the farmer and the retailer leads to inadequate access to accurate information, which negatively impacts the effectiveness of carbon reduction.
Similar to the previous sections,  under the Stackelberg and decentralized models, scenarios involving carbon trading result in significantly higher profits for the profits of both the farmer and the retailer.

\subsection{Impact of $\lambda$ and $\mu$ }
This section conducts a sensitivity analysis on the emission reduction and carbon sequestration efforts of the farmer $E_f(t)$ and the low-carbon promotion efforts of the retailer $E_r$, considering three types of cooperation models.

\begin{figure}[H]
  \centering
  \subfigure[Impact of $\lambda_f,\mu_f$ on $E_f$]{\includegraphics[width=0.49\textwidth]{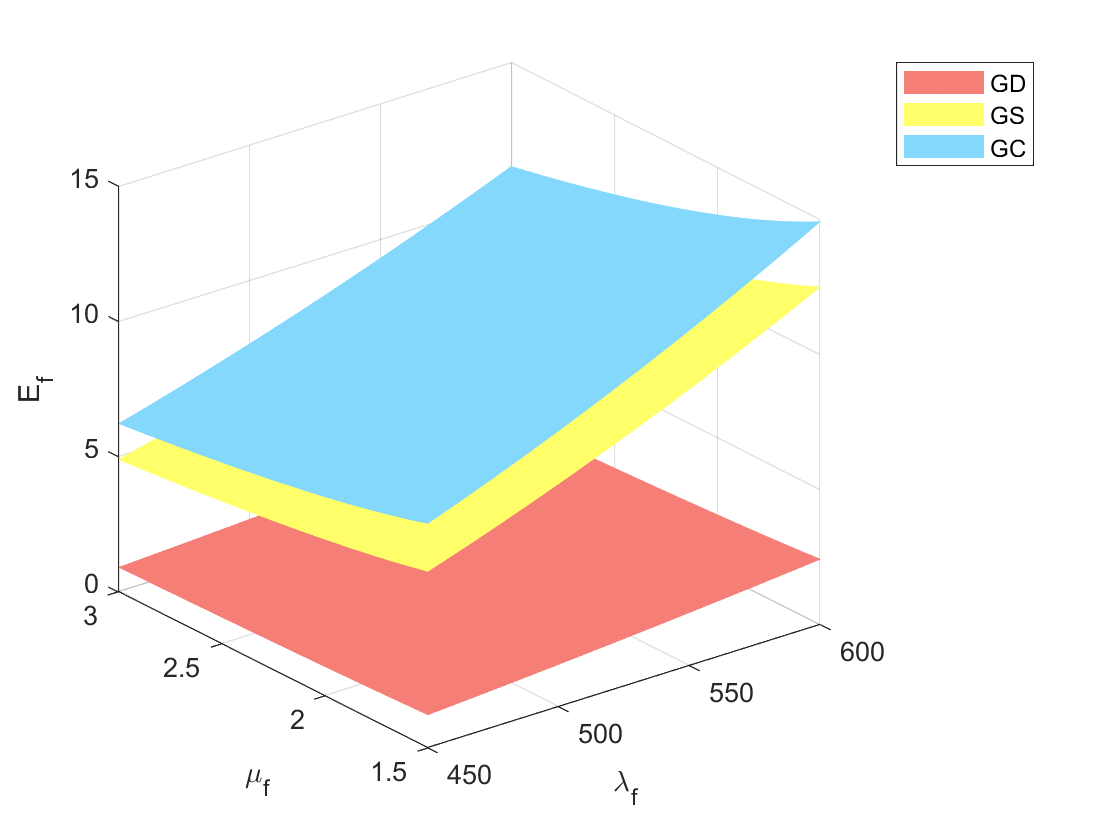}}
  \subfigure[Impact of $\lambda_r,\mu_r$ on $E_r$]{\includegraphics[width=0.49\textwidth]{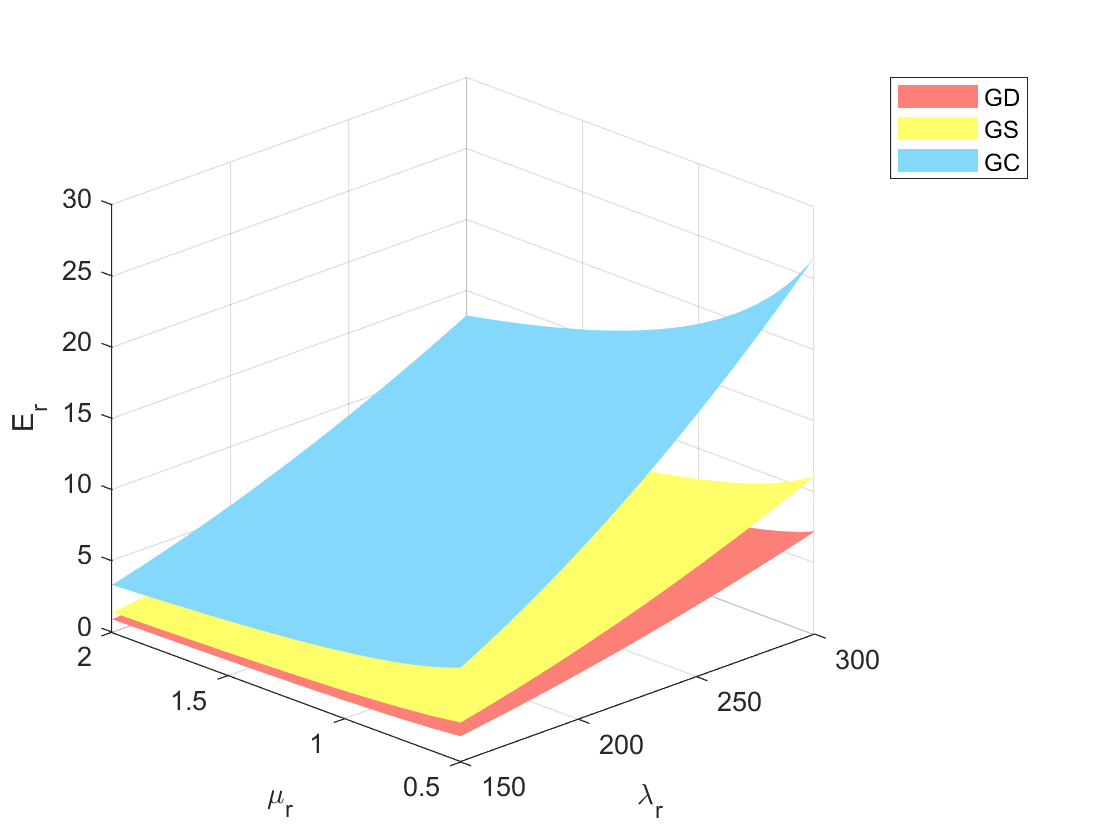}}
    \caption{Sensitivity analysis on the efforts of the farmer and the retailer}
      \label{Ef}  
\end{figure}
Fig.\ref{Ef}-a depicts the impact of $\lambda_f$ and $\mu_f$ on $E_f$, under three types of cooperation models. Fig.\ref{Ef}-a shows that as the $\lambda_f$ and $\mu_f$ increase, the emission reduction and carbon sequestration efforts of the farmer $E_f$ also increases.
Furthermore, $E_f$ is the highest in the centralized model, followed by the Stackelberg model, and the lowest in the decentralized model. 
Therefore, appropriately adjusting the relevant coefficients can incentivize both the farmer and the retailer to increase their efforts to reduce carbon emission.

Fig.\ref{Ef}-b depicts the impact of $\lambda_r$ and $\mu_r$ on $E_r$, under three types of cooperation models.
The trend of Fig.\ref{Ef}-b is similar to that of Fig.\ref{Ef}-a. 
Since the effectiveness of carbon reduction in the supply chain is significantly influenced by the efforts of the farmer and the retailer to reduce carbon emission, there is a strong correlation between the increase in these efforts and the enhanced levels of carbon reduction across three different cooperative models (as illustrated in Fig.\ref{carbon emission}-a).

\subsection{Impact of $p_c$}
This section discusses the impact of unit carbon sink price $p_c$ on the emission reduction efforts $E_f$ and profit of the farmer, considering three types of cooperation models. 

\begin{figure}[H]
  \centering
  \subfigure[Impact of carbon $p_c$ on $E_f$]{\includegraphics[width=0.49\textwidth]{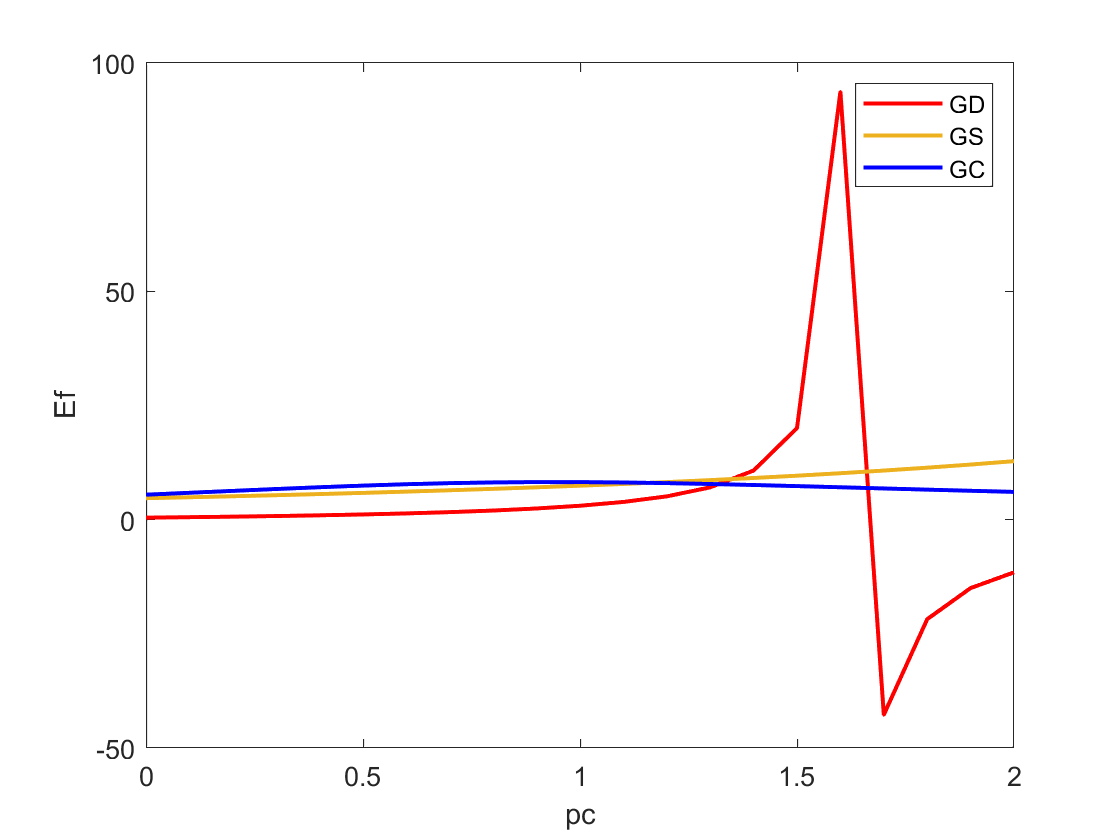}}
  \subfigure[Impact of $p_c$ on the profit of the farmer]{\includegraphics[width=0.49\textwidth]{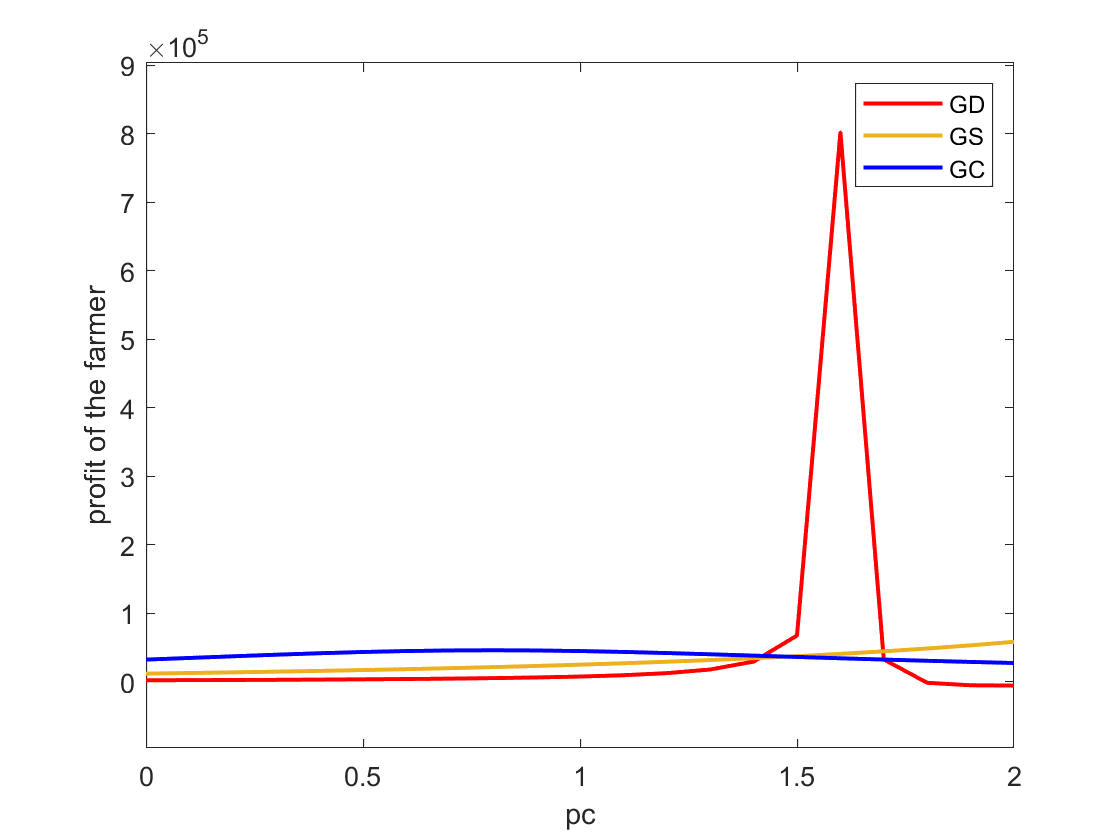}}
  \caption{Sensitivity analysis on the efforts and profit of the farmer}
  \label{pc}
\end{figure}

Fig.\ref{pc}-a depicts the impact of carbon sink price ($p_c$) on the farmer's carbon emission reduction efforts($E_f$).
In Fig.\ref{pc}-a, the most significant variations are observed under the decentralized decision model. Under this model, when $p_c$ is samll, the changes in $E_f$ are directly proportional to $p_c$. This indicates that an increase in carbon sink price can incentivize the farmer's carbon emission reduction efforts. However, as $p_c$ continues to increase and reaches a peak, further increases in $p_c$ start to inhibit the growth of $E_f$. This change can be attributed to two main reasons:
\begin{itemize}
  \item Substantial increase in carbon reduction costs due to technological limitations. And the profits of farmers also show a downward trend, when the price of carbon sinks reaches its peak(in fig.\ref{pc}-b).
  \item The farmer exhibit risk aversion behavior,  which may include concerns about a potential future decline in carbon prices, leading them to prioritize maximizing their current profits over increasing their investments in carbon emission reduction.
\end{itemize}

Fig.\ref{pc}-b depicts the impact of carbon sink price on the farmers' profit, and its trend is basically consistent with Fig.\ref{pc}-a, indicating that appropriately increasing carbon sink prices can increase farmer's profit.
\section{Conclusion}
In this paper, we explore the cooperative models between the farmer and the retailer within the context of agricultural carbon sink trading, as well as the optimal decisions on the efforts to reduce carbon emission for both parties under different cooperative models. 
In contrast to previous research that primarily focused on carbon emissions during processing and transportation in low-carbon supply chains, this paper mainly investigates carbon emissions during the production of agricultural products and innovatively incorporate carbon sink trading into the agricultural supply chain. 
Specifically, we delve into three distinct cooperative frameworks: the decentralized, the Stackelberg, and the centralized models, each accompanied by a corresponding differential game model. 
The Hamilton-Jacobi-Bellman equation is utilized to investigate the equilibrium decisions of each participant under these three cooperative models, respectively.
The main conclusions are as follows:
\begin{itemize}
    \item In the long term, the centralized model excels in all aspects, including the carbon emission reduction efforts of farmers and retailers, the carbon emission reduction level of the agricultural supply chain, and the overall profits of the supply chain, followed by the Stackelberg model, with the decentralized model showing the weakest performance. 
    This indicates that the farmer and the retailer should collaborate to reduce carbon emissions in agricultural products. 
    Specifically, in the initial stages of carbon reduction, supply chain members can adopt the Stackelberg model, where the retailers provides the farmer with essential subsidies for carbon reduction costs. 
    As the carbon emission reduction efforts progress, both participants should strengthen their cooperation and shift to a centralized cooperative model to jointly pursue deeper carbon reduction goals.
    \item Under the three types of cooperation models, performance with carbon sink trading consistently exceeds that without carbon sink trading. 
    Specifically, In the context of carbon trading, both the farmer and the retailer achieve higher profits, and the level of carbon emission reduction is also improved. 
    Additionally, an appropriately set carbon sink price significantly enhances the  efforts of the farmer to reduce carbon emission. Therefore, the government should adopt carbon sink trading as an effective incentive mechanism to promote the implementation of carbon emission reduction throughout the agricultural supply chain.
\end{itemize}

In conclusion, while this study has provided valuable insights into the effects of carbon emission reduction strategies in agricultural product supply chain, there are several avenues for further research. 
Carbon sink trading is an important policy tool in the carbon market, where there are sellers and buyers of carbon sinks. Farmers can earn additional income by selling carbon sinks, while companies in need can purchase these sinks to offset their excess carbon emissions. This study is limited to the scenario where only the farmer sells carbon sinks. Therefore, future research could further explore the involvement of companies in the carbon market and assess their potential impact on the sustainability of the agricultural product supply chain. Moreover, carbon trading policies merit consideration during the processing and transportation stages of the agricultural supply chain. Future research could elucidate the synergistic benefits of carbon trading and sink trading on supply chain efficiency and carbon emission reduction.

%
%
%
%

\newpage
\appendix
\section*{Appendix A: Proof of Proposition 1}

After discounting the profits of the farmer and the retailer over an infinite time horizon with a profit discount rate $\rho$, the objective functions for the farmer and the retailer are formulated as follows:
\setcounter{equation}{0}  
\begin{align}
\max J_{f}^{GD} & =\int_0^{\infty} e^{-\rho t}\left[p_f Q  + p_c (1+\omega E_f)Q -\frac{1}{2} \lambda_f E_f^2\right]dt, \label{qf11}\\
\max J_r^{GD} & =\int_0^{\infty} e^{-\rho t} \left[p_r D -\frac{1}{2} \lambda_r E_r^2 \right]dt.  \label{qr11}    
\end{align} 
According to (\ref{qf11}) and (\ref{qr11}), let $V_{f}^{GD}, V_{r}^{GD}$ represent the profit functions of the farmer and the retailer, respectively, then the HJB equation is as follows:
\begin{align} 
\rho V_{f}^{GD} &=\max \bigg[ p_f Q -\frac{1}{2} \lambda_f E_f^2+p_c (1+\omega E_f)Q+V_f^{{GD}^{\prime}}(\mu_f E_f+\mu_r E_r-\delta H)\bigg]dt, \label{vsD}\\
\rho V_{r}^{GD} & = \max \bigg[p_r D  -\frac{1}{2} \lambda_r E_r^2 +V_r^{{GD}^{\prime}}(\mu_f E_f+\mu_r E_r-\delta H)\bigg]dt. \label{vmD}  
\end{align} 
For ease of calculation, we define $k_1=(Q_0+ap)\theta$ and $k_2=(D_0-bp)\theta$. \\
Calculate the first-order partial derivatives with respect to $E_r$ and $E_f$ for equations (\ref{vsD}), (\ref{vmD}) respectively, and set them to zero. We can then derive the optimal efforts of the farmer and the retailer to reduce carbon emission are as follows:
\begin{align}
{E_f^{GD}}^* =\frac{p_c \theta \omega k_1 H + \mu_f V_f^{{GD}^{\prime} }}{\lambda_f}, ~~{E_r^{GD}}^* =\frac{\mu_r V_r^{{GD}^{\prime} }}{\lambda_r}. \label{em1}
\end{align}
Substituting (\ref{em1}) into (\ref{vsD}) and (\ref{vmD}), we can obtain:
\small{
\begin{equation}\label{vsd1}
\begin{aligned}
\rho V_f^{GD} &= {\frac{(p_c \omega k_1)}{2 \lambda_f}}^2 H^2+\bigg(p_f k_1+p_c k_1- \delta V_f^{{GD}^{\prime}}+\frac{\mu_f p_c \omega k_1 V_f^{{GD}^{\prime}}}{\lambda_f}\bigg)H \\
&+ {\frac{\mu_f^2 V_f^{{GD}^{\prime}}}{2 \lambda_f}}^2  +\frac{\mu_r^2 V_r^{{GD}^{\prime} }V_f^{{GD}^{\prime}}}{\lambda_r}, 
\end{aligned}
\end{equation}}
\small{
\begin{equation}
\begin{aligned}
\rho V_r^{GD} &= \bigg(p_r k_2- \delta V_r^{{GD}^{\prime}}+\frac{\mu_f p_c \omega k_1 V_r^{{GD}^{\prime}}}{\lambda_f}\bigg)H- {\frac{\mu_r^2 V_r^{{GD}^{\prime}}}{2 \lambda_r}}^2 +\frac{\mu_f^2 V_f^{{GD}^{\prime} }V_r^{{GD}^{\prime}} }{\lambda_f}. \label{vmd1}
\end{aligned}
\end{equation}}
Observing the structure of equations (\ref{vsd1}) and (\ref{vmd1}), we infer that the farmer's HJB equation is quadratic linear in $H$, while the retailer's is linear in $H$, further deriving:
\small{\begin{align}
V_f^{GD}&= A^{GD} H^2+B^{GD}H+C^{GD}, ~~
V_f^{{GD}^{\prime}}=2A^{GD}H+B^{GD}. \label{ds1}\\ 
V_r^{GD}&= M^{GD}H + N^{GD}, ~~ V_r^{{GD}^{\prime} }=M^{GD}. \label{ds}
\end{align}}
By using the method of undetermined coefficients, we can obtain:
\begin{equation*}
  \begin{aligned}  
A^{GD}&= \frac{2\rho \lambda_f+4\lambda_f \delta-4\mu_f p_c \omega k_1- \sqrt{ \triangle^{GD} }}{8 \mu_f^2}, \\
M^{GD}&=\frac{p_{m} k_{2} \lambda_{s}}{p_{c} \omega k_{1}\mu_{s}+2A^{GD}\mu_{s}^{2}+\lambda_{s}\rho-\delta\lambda_{s}}, \\
B^{GD}&=-\frac{\left(2 \lambda_r^{2} A^{GD}M^{GD}+k_1 \lambda_r \mathit{p_c}+k_1\lambda_r\mathit{p_f}\right) \lambda_f}{\lambda_r \left(\mathit{p_c} \omega k_1 \mathit{us}+2 A \, \mathit{us}^{2}-\lambda_f \rho-\delta \lambda_f\right)}, \\
C^{GD} &=\frac{B^{GD} \left(\lambda_{s}^{2}\lambda_{m} B^{GD}+2M^{GD}\lambda_{m}^{2} \lambda_{s}\right)}{2 \lambda_{m}\lambda_{s}\rho}, \\
N^{GD}&=\frac{M^{GD} \left(2 \lambda_{m}\mu_{s}^{2} B^{GD} +M^{GD}\mu_{m}^{2}\lambda_{s}\right)}{2 \lambda_{m}\lambda_{s}\rho}, 
  \end{aligned}
\end{equation*}
where $\triangle^{GD} = (4 \mu_f  p_c \omega k_1-4\lambda_f \delta-\rho\lambda_f)^2- (4\mu_f p_c\omega k_1 )^2 > 0$. \\
Furthermore, by substituting (\ref{em1}) into $\dot{H}(t)=\mu_f E_f(t)+\mu_r E_r(t)+ \mu_r E_r(t)-\delta H(t)$, we can derive the equilibrium level of carbon emission reduction when setting $\dot{H}(t)=0$, thus\\
\begin{align}
 H_d^{GD} &=\frac{\lambda_r\mu_f B^{GD}+\lambda_f\mu_r M^{GD}}{\lambda_f\lambda_r \delta-2A^{GD}\mu_s\lambda_r-p_c\omega k_1 \mu_f\lambda_r}. \nonumber  
\end{align}
The proposition 1 is proved.

\section*{Appendix A: Proof of Proposition 2}
After discounting the profits of the farmer and the retailer over an infinite time horizon with a discount rate $\rho$, the objective functions for farmer and the retailer are formulated as follows:
\setcounter{equation}{0}  
    \begin{align}
\max J_r^{GS}&=\int_0^{\infty} e^{-\rho t} \bigg[p_r D-\frac{1}{2} \lambda_r E_r^2-\frac{1}{2}x_f\lambda_f E_f^2 \bigg]dt,\label{qf22}\\
\max J_{f}^{GS} & =\int_0^{\infty} e^{-\rho t}\left[p_fQ+p_c (1+\omega E_f)Q -(1-x_f)\frac{1}{2} \lambda_f E_f^2\right]dt.\label{qr22}
\end{align}
\label{jsm1}
Using the backward induction approach, we derive the optimal decisions for each party involved in the analysis. Firstly, in the second stage, assuming that the farmer receives subsidies denoted by $x_f$, and exerts optimal carbon emission reduction efforts, denoted by $E_f$. According to (\ref{qf22}), let $V_{f}^{GS}$ represent the profit functions of the farmer, then his HJB equations is as follows: 
\begin{align} 
\rho V_{f}^{GS} & =\max \bigg[p_f Q +p_c (1+\omega E_f)Q-(1-x_f)\frac{1}{2} \lambda_f E_f^2 \nonumber\\& +V_f^{{GS}^{\prime}}(\mu_f E_f+\mu_r E_r+\mu_r E_r-\delta H)\bigg]dt. \label{vgs}
\end{align} 
For ease of calculation, we define $k_1=(Q_0+ap)\theta$ and $k_2=(D_0-bp)\theta$. \\
Calculate the first-order partial derivatives with respect to  $E_f$ for equations (\ref{vgs}) and set it to zero. We can then derive the optimal efforts of the farmer  to reduce carbon emission is as follows:
\begin{align}
{E_f^{GS}}^* &=\frac{p_c\omega k_1 H  + \mu_fV_f^{{GS}^{\prime} } }{\lambda_f}. \label{egs}
\end{align}
Next, according to (\ref{qr22}), let $V_{r}^{GS}$ represent the profit function of the retailer, then his HJB equation is as follows:
\begin{align} 
\rho V_{r}^{GS} =&\max \bigg[p_r D -\frac{1}{2} \lambda_r E_r^2-\frac{1}{2}x_f\lambda_f E_f^2+V_r^{{GS}^{\prime}}(\mu_f E_f+\mu_r E_r-\delta H)\bigg]dt. \label{vgm} 
\end{align} 
Substituting (\ref{egs}) into (\ref{vgm}), we then find the first-order partial derivatives with respect to $E_r$ and $x_f$ and set them to zero. From this, we can obtain the optimal decision of the retailer in the first stage, as shown below:
\begin{align}
{E_r^{GS}}^*=\frac{\mu_r V_r^{{GS}^{\prime}}}{\lambda_r}, ~~x_f^{GS^{*}}=\frac{2V_r^{GS^{\prime}}-V_f^{GS^{\prime}}-p_c\omega k_1 H/\mu_f}{2V_r^{GS^{\prime}}+V_f^{GS^{\prime}}+p_c\omega k_1 H/\mu_f}. \label{xs}
\end{align}
Substitute (\ref{egs}), (\ref{xs}) into (\ref{vgs}) and (\ref{vgm}) we can get:
\begin{align}
\rho V_f^{GS} &= \frac{\mu_f (p_c\omega k_1)^2}{4\lambda_f}H^2+\bigg(p_f k_1+p_c k_1-\delta V_f^{{GD}^{\prime}}+\frac{\mu_f p_c\omega k_1 (V_f^{{GD}^{\prime}}+V_r^{{GD}^{\prime}})}{2\lambda_f}\bigg)H\nonumber\\
&+\frac{\mu_r^2 V_f^{{GD}^{\prime}}V_r^{{GD}^{\prime}}}{\lambda_r}
+\frac{\mu_f^2(2V_r^{{GD}^{\prime}}V_f^{{GD}^{\prime}}+{V_f^{{GD}^{\prime}}}^{2})}{4\lambda_f}, \label{vgs1}
\end{align}
\begin{align}
\rho V_r^{GS} &=\frac{\mu_f (p_c\omega k_1)^2}{8\lambda_f}H^2+\bigg(p_r k_2+\frac{\mu_f p_c\omega k_1 ( V_f^{{GD}^{\prime}}+2V_r^{{GD}^{\prime}})}{4\lambda_f}-\delta V_r^{{GD}^{\prime}}\bigg)H\nonumber\\&+\frac{\mu_r^2 {V_r^{{GD}^{\prime}}}^{2}}{2\lambda_r}+\frac{\mu_f^2 (V_f^{{GD}^{\prime}}+2V_r^{{GD}^{\prime}})^2}{8\lambda_f}. \label{vgm1}
\end{align}
Observing the structures of (\ref{vgs1}) and (\ref{vgm1}), we can see that the profit functions of the farmer and the retailer are quadratic linear functions of $H$, as shown below:
\begin{align}
V_f^{GS}&= A^{GS} H^2+B^{GS}H +C^{GS}, ~~V_f^{{GS}^{\prime}}=2A^{GS}H+B^{GS}. \label{vgs2}\\
V_r^{GS}&= M^{GS} H^2+N^{GS}H+T^{GS}, ~~V_r^{{GS}^{\prime}}=2M^{GS}H+N^{GS}. \label{vgm2}
\end{align}
By using the method of undetermined coefficients, we can obtain:
\begin{equation*}
\begin{aligned}
A^{GS}&=\frac{2\lambda_r\lambda_f\rho \delta-p_c\omega\ k_1\lambda_r\rho\mu_f-\lambda_r\lambda_f\rho^3-2\mu_f^2\lambda_r M^{GS}-4\mu_r^2\lambda_f M^{GS}-\sqrt{\triangle^{GS1}}}{2\mu_f^2\lambda_r},  \\
M^{GS}&=\frac{2\lambda_f\lambda_r\rho\delta-p_c\omega k_1\rho\lambda_r\mu_f-\lambda_f\lambda_r\rho^3-2A^{GS}\lambda_r\mu_f^2-\sqrt{\triangle^{GS2}}}{4(\lambda_r\mu_f^2+\lambda_f\mu_r^2)}, \\
B^{GS}&=\frac{\lambda_f\lambda_r((p_c +p_f)k_1+p_r k_2)}{\lambda_f\lambda_r(\rho-\delta)-p_c\omega\lambda_r\mu_f k_1-2A^{GC}(\lambda_r\mu_f^2+\lambda_f\mu_r^2)}, \\
C^{GS}&=\frac{{B^{GS}}^2\mu_f^2\lambda_r+B^{GS}N^{GS}(4\mu_r^2\lambda_f+2\mu_f^2\lambda_r)}{4\rho^3\lambda_f\lambda_r}, \\
N^{GS}&=\frac{2\lambda_r\mu_f^2B^{GS}(A^{GS}+2M^{GS}-4\lambda_r\lambda_f p_r\rho^2 k_2-p_c\omega k_1\lambda_r\mu_f\rho B^{GS})}{4\rho^3\lambda_f\lambda_r-2p_c\omega k_1\mu_f\rho\lambda_r+4\lambda_f\lambda_r\delta\rho+4\mu_f^2\lambda_r(A^{GS}+2M^{GS})+8\mu_r^2\lambda_f M^{GS}}, \\
F^{GS}&=\frac{\mu_f^2\lambda_r({B^{GS}}^2+ 4N^{GS}B^{GS})+4{N^{GS}}^2(\lambda_r\mu_f^2+\lambda_f\mu_r^2)}{8\rho^3\lambda_r
\lambda_f}, 
\end{aligned}
\end{equation*}
where $k_1=(Q_0+ap)\theta$, $k_2=(D_0-bp)\theta$,\\
\begin{equation}
    \begin{aligned}
\triangle^{GS1}&=2(\mu_f \omega k_1 p_c\rho)^{2} +2k_1 p_c\omega\lambda_f\mu_f \rho^{2}(\rho^{2}-2\delta) +\lambda_f\rho^{3}(\lambda_f\rho^{3}+4M\mu_f^{2}-4\lambda_f\rho\delta)\\
&+4M^{2}\mu_f^{4}-8M\lambda_f\rho\delta\mu_f^{2}+4\lambda_f^{2}+4(\lambda_f\rho\delta)^{2}+8M\mu_r^{2}\lambda_f(\lambda_r\rho\omega p_c k_1+\lambda_f\lambda_r\rho^{3}\\
&+2M\mu_f^{2}\lambda_r+2M\mu_r^{2}\lambda_f-2\lambda_r\lambda_f\rho\delta)>0, and\\
\triangle^{GS2}&=2\lambda_r(k_2 p_c\rho\omega)^{2}(\lambda_r\mu_f^{2}+\lambda_f)+2k_2\lambda_r\lambda_f p_c \rho\omega(\lambda_r\rho^{3}\mu_f-2A\mu_r^{2}-2\lambda_r\rho\delta)\\
&+4\lambda_r\lambda_f\rho^{2}(\rho^{3}+\rho^{2}\delta+\delta^{2})+4A\lambda_r\lambda_f\mu_f^{2}(\lambda_r\rho^{3}-A\mu_r^{2}-2\lambda_r\rho\delta)>0.\nonumber
    \end{aligned}
\end{equation}
Furthermore, by substituting (\ref{egs}) into $\dot{H}(t)=\mu_f E_f(t)+\mu_r E_r(t)+ \mu_r E_r(t)-\delta H(t)$, we can derive the equilibrium level of carbon emission reduction when setting $\dot{H}(t)=0$, thus\\
\begin{align}
H_d^{GS}&=\frac{2(\lambda_f\rho\mu_r +\lambda_r\mu_f)N^{GS}+\lambda_r\mu_f B^{GS}}{4\lambda_r\mu_f M^{GS}-\rho \lambda_r p_c \omega k_1-4\lambda_f\mu_r\rho M^{GS}-2\lambda_r\mu_f A^{GS}+2\lambda_r\lambda_f\rho\delta}.\nonumber
\end{align}
The proposition 2 is proved. 
\section*{Appendix A: Proof of Proposition 3}
\setcounter{equation}{0}  
After discounting the profits of the farmer and the retailer over an infinite time horizon with a discount rate $\rho$, the overall objective functions for farmer and the retailer are formulated as follows:
    \begin{align}
\max J_{fr}^{GC}& =\int_0^{\infty} e^{-\rho t}\bigg[p_f Q+p_r D + p_c (1+\omega E_f)Q-\frac{1}{2}\lambda_f E_f^2 -\frac{1}{2}\lambda_r E_r^2\bigg] dt. \label{qfrr}
 \end{align}
 According to (\ref{qfrr}), let $V_{fr}^{GC}$ represent the profit function of the agricultural product supply chain. The corresponding HJB equation is as follows:
\begin{align}
\rho V_{fr}^{GC}=\max \bigg[p_fk_1 H +&p_r k_2 H +p_c (1+\omega E_f)Q-\frac{1}{2} \lambda_f E_f^2\nonumber \\
&-\frac{1}{2} \lambda_r E_r^2+V_{fr}^{{GC}^{\prime}}(\mu_f E_f+\mu_r E_r-\delta H)\bigg]dt, \label{vsmc}        
\end{align}
For ease of calculation, we define $k_1=(Q_0+ap)\theta$ and $k_2=(D_0-bp)\theta$. \\
Calculate the first-order partial derivatives with respect to $E_r$ and $E_f$ for equations (\ref{vsmc}) respectively, and set them to zero. We can then derive the optimal efforts of the farmer and the retailer to reduce carbon emission are as follows:
\begin{align}
{E_f^{GC}}^* =\frac{\mu_fV_{fr}^{{GC}^{\prime} }+p_c\omega k_1 H}{\lambda s}, ~~{E_r^{GC}}^* =\frac{ \mu_r V_{fr}^{{GC}^{\prime} }}{\lambda_r}. \label{emC}
\end{align}
Substituting  (\ref{emC}) into (\ref{vsmc}), we can get:
\begin{align}
\rho V_{fr}^{GC} & = {\frac{(p_c \omega k_1)}{2 \lambda_f}}^2 H^2+\bigg[(p_f+p_c)k_1 +p_rk_2 +\frac{p_c \omega k_1 \mu_f V_{fr}^{{GC}^{\prime} }}{\lambda_f}-\delta V_{fr}^{{GC}^{\prime} }\bigg]H\nonumber  \\
& +V_{fr}^{{GC}^{\prime}} \bigg( \frac{\mu_f^2 V_{fr}^{{GC}^{\prime} }}{2\lambda_f }+\frac{\mu_r^2 V_{fr}^{{GC}^{\prime} }}{2\lambda_r}\bigg). \label{vsmcd2}  
\end{align}
Observing the structure of (\ref{vsmcd2}), it can be inferred that the profit function of the system is a quadratic linear function with respect to $H$, as shown below:
\begin{align}
V_{fr}^{GC}= A^{GC}H^2 + B^{GC} H+C^{GC}, 
~~V_{fr}^{{GC}{\prime}}=2A^{GC} H+B^{GC}. \label{vsmcd3}
\end{align}
By using the method of undetermined coefficients, we can obtain:
\begin{equation*}
    \begin{aligned}
 A^{GC}&=\frac{2\delta\lambda_f-\rho\lambda_f-2\lambda_r\mu_fp_c\omega k_1-\sqrt{\triangle^{GC}}}{4\lambda_r \mu_f^2+4\lambda_f \mu_r^2}, \\
B^{GC}&=\frac{\lambda_f\lambda_r((p_c +p_f)k_1+p_r k_2)}{\lambda_f\lambda_r(\rho-\delta)-p_c\omega\lambda_r\mu_f k_1-2A^{GC}(\lambda_r\mu_f^2+\lambda_f\mu_r^2)}, \nonumber
\\
C^{GC}&=\frac{B^{{GC}^{2}}(\lambda_r\mu_f^2+\lambda_f\mu_r^2)}{2\lambda_r\lambda_f\rho},    
    \end{aligned}
\end{equation*}
where $\triangle^{GC}=4\lambda_r(p_c \omega k_1)^2(2\lambda_r\mu_f^2+\lambda_f\mu_r^2+\lambda_r^2\lambda_f(\rho-2\delta)(4p_c\omega k_1\mu_f+\rho-2\delta))>0$. \\
Furthermore, by substituting (\ref{emC}) into $\dot{H}(t)=\mu_f E_f(t)+\mu_r E_r(t)+ \mu_r E_r(t)-\delta H(t)$, we can derive the equilibrium level of carbon emission reduction when setting $\dot{H}(t)=0$, thus\\
\begin{align}
H_d^{GC} &=\frac{(\lambda_f\mu_r+\lambda_r\mu_f)B^{GC}}{\lambda_f\lambda_r\delta-\lambda_r p_c\omega k_1-2(\lambda_f\mu_r+\lambda_r\mu_f)A^{GC}}.\nonumber
\end{align}
The proposition 3 is proved.
\end{document}